\documentclass[11pt]{article}
\usepackage[utf8]{inputenc}
\usepackage{amsthm, amsmath, amssymb, amsfonts, url, booktabs, tikz, setspace, fancyhdr, bm}
\usepackage{hyperref}
\usepackage{geometry}
\usepackage{hyperref, enumerate}
\usepackage[shortlabels]{enumitem}
\usepackage{enumitem}
\usepackage[english]{babel}
\usepackage[capitalise]{cleveref}
\usepackage{bbm,tkz-graph,subcaption}
\usepackage{csquotes}
\usepackage{mathrsfs}
\usepackage{mathabx}
\usetikzlibrary{patterns}
\usetikzlibrary{shapes}
\usepackage{bbm,wrapfig}

\usepackage{algorithm}
\usepackage{algorithmic}
\usepackage[utf8]{inputenc} 
\usepackage[T1]{fontenc}
\usepackage{amsfonts}
\usepackage{amscd}
\usepackage{graphicx}
\usepackage{enumitem}
\usepackage{array, booktabs, makecell}
\usepackage{verbatim}
\usepackage{hyperref}
\usepackage{amsmath,caption}
\usepackage{url,pdfpages,xcolor,framed,color}
\usepackage{todonotes}
\usepackage{comment}

\newtheorem{thr}{Theorem}[section]
\newtheorem{q}[thr]{Question}

\newtheorem{conj}[thr]{Conjecture}
\theoremstyle{definition}

\newtheorem*{defi*}{Definition}

\newtheorem{obs}[thr]{Observation}

\newcommand*{\myproofname}{Proof}

\tikzstyle{P} = [draw, circle, black, fill, inner sep = 0pt, minimum width = 3pt]
\tikzstyle{every loop} = []

\title{A note on $\bar{X}$-coloring and $\hat{A}$-coloring 4-regular graphs}

\date{}
\author{
Jorik Jooken\thanks{Department of Computer Science, KU Leuven Campus Kulak-Kortrijk, 8500 Kortrijk, Belgium. Jorik Jooken is supported by an FWO grant with grant number 1222524N. \protect\href{mailto:jorik.jooken@kuleuven.be}{\protect\nolinkurl{jorik.jooken@kuleuven.be}}}
}

\begin{document}
\maketitle
\begin{abstract}
\noindent Let $\partial_H(u)$ be the set of edges incident with a vertex $u$ in the graph $H$. We say that a graph $G$ is $H$-colorable if there exist total functions $f : E(G) \rightarrow E(H)$ and $g : V(G) \rightarrow V(H)$ such that $f$ is a proper edge-coloring of $G$ and for each vertex $u \in V(G)$ we have  $f(\partial_G(u))=\partial_H(g(u))$. Let $\bar{X}$ be the graph obtained by adding three parallel edges between two degree one vertices of the graph $K_{1,4}$. Let $\hat{A}$ be the graph obtained by adding two pendant edges to two different vertices of a triangle and then adding two edges between the degree two vertex and the two adjacent degree three vertices. Malnegro and Ozeki [\emph{Discrete Math.} \textbf{347}(3):113844 (2024)] asked whether every 4-regular graph with an even number of vertices and an even cycle decomposition of size 3 admits an $\bar{X}$-coloring or an $\hat{A}$-coloring and whether every 2-connected planar 4-regular graph with an even number of vertices admits such a coloring. Additionally, they conjectured that for every 2-edge-connected simple cubic graph $G$ with an even number of edges, the line graph $L(G)$ is $\bar{X}$-colorable. In this short note, we discuss two algorithms for deciding whether a graph $G$ is $H$-colorable. We give a negative answer to the two questions and disprove the conjecture by finding suitable graphs, as verified by two independent algorithms.


\end{abstract}

\section{Introduction}

Among the regular graphs, cubic graphs have received the most attention in the context of $H$-colorings. The most notable example here is the Petersen Coloring Conjecture due to Jaeger~\cite{J88}.

\begin{conj}[\cite{J88}]
Any 2-edge-connected cubic graph is $P$-colorable, where $P$ denotes the Petersen graph.
\end{conj}

This conjecture is very important, because its veracity would imply a number of important conjectures such as the Berge-Fulkerson conjecture~\cite{F71} and the Cycle Double Cover conjecture~\cite{S79,S73}.

\begin{figure}[h!]
  \begin{subfigure}[b]{0.35\textwidth}
    \renewcommand{\arraystretch}{0.9}
    \centering
    \includegraphics[width=0.6\linewidth]{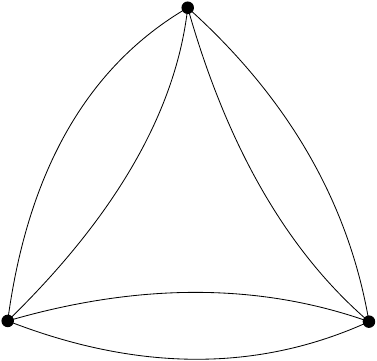}
    \caption{}
    \label{labelHere}
  \end{subfigure}%
  \begin{subfigure}[b]{0.25\textwidth}
    \renewcommand{\arraystretch}{0.5}
    \centering
    \includegraphics[width=0.6\linewidth]{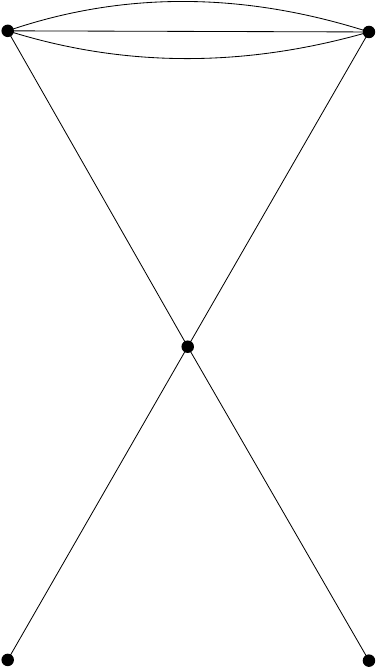}
    \caption{}
    \label{labelHere}
  \end{subfigure}%
  \begin{subfigure}[b]{0.35\textwidth}
    \renewcommand{\arraystretch}{0.9}
    \centering
    \includegraphics[width=0.6\linewidth]{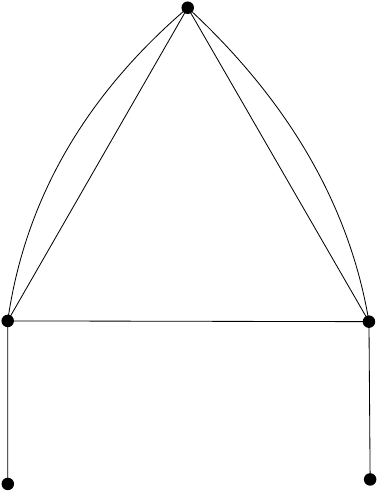}
    \caption{}
    \label{labelHere}
  \end{subfigure}%
  \caption{The graphs $2C_3$ (a), $\bar{X}$ (b) and $\hat{A}$ (c).}
\label{fig:threeGraphs}
\end{figure}

$H$-colorings of $k$-regular graphs generalize the classical notion of proper edge-colorings of $k$-regular graphs, since a $k$-regular graph $G$ has a proper $k$-edge-coloring if and only if $G$ admits a $K_{1,k}$-coloring. The 4-regular case has recently been studied by Malnegro and Ozeki~\cite{MO24}. They identified a number of small graphs $H$ for which it is interesting to ask whether a given 4-regular graph is $H$-colorable and relate being $H$-colorable with well-known properties such as having an even 2-factor and having a proper 4-edge-coloring. Three such graphs that they identified are $2C_3$, $\bar{X}$ and $\hat{A}$ (shown in~\cref{fig:threeGraphs}).

They showed that a 4-regular graph $G$ is $2C_3$-colorable if and only if $G$ has an even cycle decomposition of size 3. For $\bar{X}$ and $\hat{A}$ no such equivalences are known, but they showed that (trivially) every 4-edge-colorable graph is also $\bar{X}$-colorable and $\hat{A}$-colorable and that every $\hat{A}$-colorable graph is also $2C_3$-colorable. They asked the following questions:
\begin{q}[\cite{MO24}]\label{ques:planar}
Does every 2-connected planar 4-regular graph with an even number of vertices admit an $\bar{X}$-coloring or an $\hat{A}$-coloring?
\end{q}
\begin{q}[\cite{MO24}]\label{ques:cycleDecomposition}
Does every 4-regular graph with an even number of vertices and an even cycle decomposition of size 3 admit an $\bar{X}$-coloring or an $\hat{A}$-coloring?
\end{q}
\noindent and made the following conjecture:
\begin{conj}[Conjecture 21 in~\cite{MO24}]\label{conj:lineGraphXBar}
If $G$ is a 2-edge-connected simple cubic graph with an even number of edges, then the line graph $L(G)$ of $G$ admits an $\bar{X}$-coloring.
\end{conj}

In the remainder of this note, we consider graphs that might have multiple edges, but no loops. In~\cref{subsec:notationAndDefinitions} we summarize the notation and definitions that will be used. We discuss an algorithm for deciding whether a graph $G$ is $H$-colorable in~\cref{sec:algo}. In~\cref{sec:mainResult} we discuss the main result of this note. We find suitable graphs that provide a negative answer for Questions~\ref{ques:planar} and~\ref{ques:cycleDecomposition} and a counterexample to~\cref{conj:lineGraphXBar} and verify this using the algorithm from~\cref{sec:algo} and an independent verification algorithm. Finally in~\cref{sec:futureResearch} we discuss an interesting avenue for future research.

\subsection{Notation and definitions}\label{subsec:notationAndDefinitions}
A \textit{partial function} is a function that is defined for a (possibly empty) subset of elements from its domain. A \textit{total function} is a partial function that is defined for all elements of its domain. For an integer $k \geq 1$ we denote by $[k]$ the set $\{1,2,\ldots,k\}$. A \textit{proper $k$-edge-coloring} of a graph $G$ is a total function $f : E(G) \rightarrow [k]$ (in fact the co-domain can be any $k$-element set) such that for adjacent edges $e_1, e_2 \in E(G)$ we have $f(e_1) \neq f(e_2)$. A graph which admits a proper $k$-edge-coloring is said to be \textit{$k$-edge-colorable}. The \textit{chromatic index} of a graph is the smallest integer $k$ such that the graph is $k$-edge-colorable. The set of edges incident with a vertex $u$ in the graph $H$ is denoted by $\partial_H(u)$. We say that a graph $G$ is \textit{$H$-colorable} if there exist total functions $f : E(G) \rightarrow E(H)$ and $g : V(G) \rightarrow V(H)$ such that $f$ is a proper edge-coloring of $G$ and for each vertex $u \in V(G)$ we have $\{f(e)~|~e \in \partial_G(u)\}=\partial_H(g(u))$ (we lift the definition of $f$ in the usual way and abbreviate this as $f(\partial_G(u))=\partial_H(g(u))$). A \textit{cycle} is a connected graph in which every vertex has degree two. An \textit{even cycle decomposition} of a graph $G$ is a partition of $E(G)$ into edge-disjoint cycles of even length. The \textit{size} of an even cycle decomposition is the minimum number of colors needed to color all cycles in the even cycle decomposition such that no two cycles with the same color share a vertex.

\section{An algorithm for deciding if a graph $G$ is $H$-colorable}\label{sec:algo}
In this section we describe a relatively straightforward recursive algorithm that can decide whether a graph $G$ is $H$-colorable. Given graphs $G$ and $H$ and a partial function $f : E(G) \rightarrow E(H)$, we say that a color $c \in E(H)$ is \textit{compatible} with an edge $u_1u_2 \in E(G)$ for which $f$ is not defined if for each vertex $v \in \{u_1, u_2\}$, $c$ is different from $f(vv')$ for all $vv' \in E(G)$ for which $f$ is defined and there exists a vertex $w \in V(H)$ such that the set $\{c\} \cup \{f(vv')~|~vv' \in E(G)\text{ and f is defined for }vv'\}$ is a subset of $\partial_H(w).$ The recursive algorithm is called with two partial functions $f : E(G) \rightarrow E(H)$ and $g : V(G) \rightarrow V(H)$ that are not defined for any of its domain elements. In each recursion step, the algorithm chooses an edge $u_1u_2 \in E(G)$ for which $f$ is not defined and $u_1u_2$ has the least number of compatible colors. The algorithm assigns a compatible color $c \in E(H)$ to this edge in all possible ways and extends $f$ by defining $f(u_1u_2)=c.$ If after doing this, $f$ is defined for all edges incident with $u_1$ (or $u_2$), the algorithm extends $g$ by defining $g(u_1)=w$ for some vertex $w \in V(H)$ such that $f(\partial_G(u_1))=\partial_H(w)$ (and similarly for $u_2$). This simple recursive step is repeated until both $f$ and $g$ are total functions (in which case $G$ is $H$-colorable) or until there are no more possibilities to extend $f$ and $g$ (in which case $G$ is not $H$-colorable). The rationale for choosing the next edge to color as the edge with the least number of compatible colors follows the fail-first principle and some preliminary testing revealed that this was indeed a good choice. The pseudo code of the recursive algorithm is shown in~\cref{algo:decisionAlgo}.

\begin{algorithm}[ht!]
\caption{isColorable(Graph $G$, Graph $H$, Partial function $f : E(G) \rightarrow E(H)$, Partial function $g : V(G) \rightarrow V(H)$)}
\label{algo:decisionAlgo}
  \begin{algorithmic}[1]
		\STATE // Returns TRUE if and only if the partial functions $f$ and $g$ can be extended such that they form an $H$-coloring of $G.$

        \IF{$f$ and $g$ are total functions}
            \RETURN TRUE
        \ENDIF
        \STATE Let $u_1u_2 \in E(G)$ be an edge for which $f$ is not defined that has the least number of compatible colors
        \FOR{each color $c \in E(H)$ compatible with $u_1u_2$}
                \STATE $f' \gets extend(f,u_1u_2,c)$
                \STATE $g' \gets g$
                \FOR{$v \in \{u_1, u_2\}$}
                    \IF{$f'$ is defined for all edges incident with $v$}
                        \STATE Let $w \in V(H)$ be a vertex for which $f(\partial_G(v))=\partial_H(w)$
        			\STATE $g' \gets extend(g',v,w)$
                    \ENDIF
                \ENDFOR
                \IF{$isColorable(G,H,f',g')$}
                    \RETURN TRUE
                \ENDIF
        \ENDFOR
		\RETURN FALSE
  \end{algorithmic}
\end{algorithm}

\section{Negative answers for Questions~\ref{ques:planar} and~\ref{ques:cycleDecomposition} and a counterexample to~\cref{conj:lineGraphXBar}}\label{sec:mainResult}

\begin{figure}[h!]
	\centering
  \includegraphics[width=0.6\linewidth]{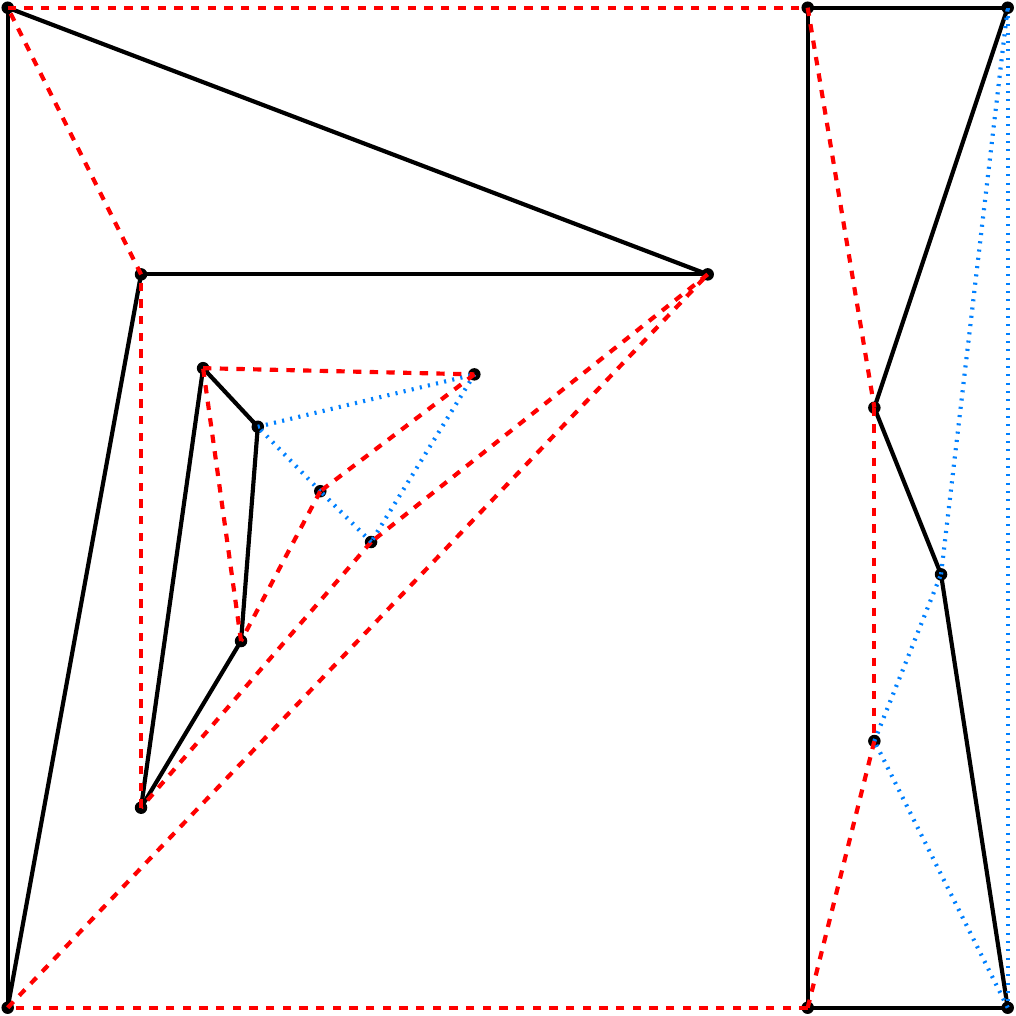}
\caption{A 2-connected planar 4-regular graph with an even number of vertices and an even cycle decomposition of size 3 that is not $\bar{X}$-colorable nor $\hat{A}$-colorable.}
\label{fig:planarGraph}  
\end{figure}

We implemented~\cref{algo:decisionAlgo} to decide if a given graph $G$ is $H$-colorable. We used the planar graph generator \textit{plantri}~\cite{BM07} and the aforementioned algorithm to generate all simple 2-connected planar 4-regular graphs with an even number of vertices that are $2C_3$-colorable, i.e. have an even cycle decomposition of size 3, on a fixed number of vertices. Since all 4-regular graphs which are 4-edge-colorable are also $\bar{X}$-colorable and $\hat{A}$-colorable, we further restricted the search to graphs with chromatic index 5 using the \textit{nauty} package~\cite{MP14}. Note that the chromatic index of a simple 4-regular graph is at most 5, because of Vizing's theorem~\cite{V65}. Using~\cref{algo:decisionAlgo}, we found that the graph shown in~\cref{fig:planarGraph} is the smallest simple 2-connected planar 4-regular graph with an even number of vertices that is not $\bar{X}$-colorable nor $\hat{A}$-colorable (the graph has 18 vertices). Moreover, the graph also has an even cycle decomposition of size 3 (indicated by the 3 different edge styles in this figure). Hence, this graph simultaneously provides a negative answer for Questions~\ref{ques:planar} and~\ref{ques:cycleDecomposition}. We made this graph available on House of Graphs~\cite{HOG} at \url{https://houseofgraphs.org/graphs/50558}.

\begin{figure}[h!]
  \begin{subfigure}[b]{0.5\textwidth}
    \renewcommand{\arraystretch}{0.9}
    \centering
    \includegraphics[width=0.6\linewidth]{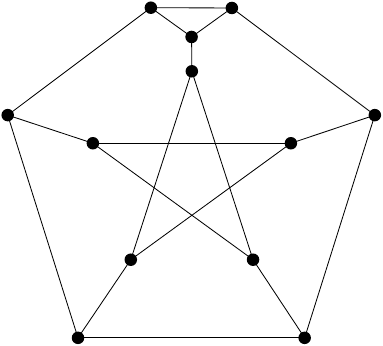}
    \caption{}
    \label{labelHere}
  \end{subfigure}%
  \begin{subfigure}[b]{0.5\textwidth}
    \renewcommand{\arraystretch}{0.9}
    \centering
    \includegraphics[width=0.6\linewidth]{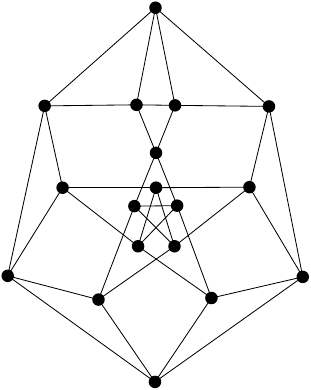}
    \caption{}
    \label{labelHere}
  \end{subfigure}%
  \caption{The Tietze graph (a) and its line graph (b).}
\label{fig:TietzeAndLineGraph}
\end{figure}

We also used the cubic graph generator \textit{snarkhunter}~\cite{BGHM13} and the \textit{nauty} package~\cite{MP14} to generate until order 12 all 2-edge-connected simple cubic graphs with an even number of edges. This enabled us to show that the smallest such graph whose line graph does not admit an $\bar{X}$-coloring is the Tietze graph, which can be obtained by replacing a vertex in the Petersen graph by a triangle. This graph (shown in~\cref{fig:TietzeAndLineGraph}) is a counterexample to~\cref{conj:lineGraphXBar} and made available at \url{https://houseofgraphs.org/graphs/1368}.

\subsection{Independent verification}
Since we conclude that the graphs are not $H$-colorable based on the output of an algorithm, it is very important to ensure that the algorithm is implemented correctly. To this end we also implemented another algorithm that can decide whether a graph $G$ is $H$-colorable. This second algorithm is a recursive algorithm that starts with two partial functions $f : E(G) \rightarrow E(H)$ and $g : V(G) \rightarrow V(H)$ that are not defined for any of its domain elements. In each recursion step, the algorithm extends $g$ by choosing an arbitrary vertex $v \in V(G)$ for which $g$ is not yet defined and tries all valid assignments $g(v)=w$ that are consistent with the current coloring. After extending $g$, the function $f$ is also extended in all possible ways for all edges $vv' \in E(G)$ such that $f$ was not previously defined for $vv'$. This second algorithm tends to be slower than~\cref{algo:decisionAlgo}, but it is useful to obtain an independent verification. Using this second algorithm, we repeated the same methodology as described in~\cref{sec:mainResult} and we obtained the same outcome for all graphs with both algorithms. We make the source code of both algorithms publicly available at~\url{https://github.com/JorikJooken/HColoring}.

\section{Future research}\label{sec:futureResearch}
In this note, we gave a negative answer to two questions and disproved a conjecture by Malnegro and Ozeki~\cite{MO24}. This narrows down the avenues that should be considered for future research. One such avenue is the analogue conjecture of~\cref{conj:lineGraphXBar}, where $\bar{X}$ is replaced by $\hat{A}$:

\begin{conj}[Conjecture 16 in~\cite{MO24}]\label{conj:lineGraphAHat}
If $G$ is a 2-edge-connected simple cubic graph with an even number of edges, then the line graph $L(G)$ of $G$ admits an $\hat{A}$-coloring.
\end{conj}
Malnegro and Ozeki gave partial solutions for this conjecture by showing that the line graphs of a large number of snarks (cubic graphs which are not 3-edge-colorable and satisfy additional constraints on the girth and edge-connectivity) are $\hat{A}$-colorable. Using both~\cref{algo:decisionAlgo}, the independent verification algorithm and the cubic graph generator \textit{snarkhunter}~\cite{BGHM13}, we verified for all simple 2-edge-connected cubic graphs with an even number of edges and order at most 24 that their line graph admits an $\hat{A}$-coloring. This leads to the following observation:
\begin{obs}
If $G$ is a 2-edge-connected simple cubic graph with an even number of edges, such that the line graph $L(G)$ of $G$ does not admit an $\hat{A}$-coloring, then $L(G)$ has order at least 42.
\end{obs}
This provides further evidence towards~\cref{conj:lineGraphAHat} and this conjecture deserves further attention.

\section*{Acknowledgements}
The author would like to thank Jan Goedgebeur and Davide Mattiolo for useful discussions and suggestions related to the presentation of this note.
\bibliographystyle{abbrv}
\bibliography{ref}

\end{document}